\documentclass[11pt]{article}
\usepackage[centertags]{amsmath}
\usepackage{amsfonts}
\usepackage{amssymb}
\usepackage{amsthm}
\usepackage{newlfont}
\usepackage{graphicx}
\usepackage{multirow}
\usepackage{rotating}
\usepackage{epstopdf}

\newtheorem{theorem}{Theorem}

\newtheorem{lemma}{Lemma}

\newtheorem{corollary}{Corollary}
\begin{document}
\thispagestyle{empty}
%
%

\title{ Existence and Uniqueness of the Solution of Two-dimensional Fuzzy Volterra Integral Equation with Piecewise Kernel  }
\author{Samad Noeiaghdam}
  \date{}
 \maketitle

\begin{center}
\scriptsize{ Institute of Mathematics, Henan Academy of Sciences, Zhengzhou, 450046, China. snoei@hnas.ac.cn\\
}
\end{center}
\begin{abstract}
This study investigates the existence and uniqueness of solutions to Volterra integral equations with discontinuous kernels in both linear and nonlinear cases. The problem is two-dimensional, and the collocation method is employed to analyze the equations. The research aims to provide a comprehensive understanding of the solution properties of these integral equations, which are crucial in various mathematical and physical applications. By examining the existence and uniqueness of solutions, this study contributes to the development of numerical methods for solving Volterra integral equations with discontinuous kernels. The findings of this research have the potential to impact various fields, including physics, engineering, and economics, where integral equations play a significant role in modeling complex phenomena.

\vspace{.5cm}{\it keywords: Fuzzy Volterra integral equation; Existence and Uniqueness of the solution; Collocation method  }
\end{abstract}
\section{Introduction}

Engineering and physics both include several uses for VIEs \cite{21}. The load balancing issue, energy storages, and power systems with renewable and storage batteries are only a few of the issues that arise. Additional applications can be found in the basic storage properties such as capacity, efficiency, number of cycles, discharge/charge rate, load distribution between the available storages based on forecasted electric load, and generation from conventional and renewable energy sources. This particular class of integral equations was initially described by Sidorov and Lorenzi. The VIEs for renewable and diesel generating in energy storage have been covered by Sidorov et al. \cite{48}. Tao et al. have studied the use of neural networks to anticipate load and pollution in \cite{63}. Sidorov et al. in \cite{49} examined the solution of nonlinear VIE systems with a focus on a few applications to power system operation. For the best power flow, Domyshev et al. revised the two-step optimization process \cite{12}. Sidorov et al. approximated the Volterra model of energy storage to forecast the demand in power networks in \cite{57}. The precision of the loading batteries has finally been managed by Noeiaghdam et al. in \cite{37} using the VIE with discontinuous kernel. Also many other studies on the numerical solution of Volterra integral equation with discountinuous kernel can be found in \cite{s100,s101,s102,s103,s104,s105,s106,s107}.

Because of the mentioned important applications, we study the following two dimensional linear and nonlinear FVIE with piecewise kernel

\begin{equation}\label{1}
z(s,t) = g(s,t) \oplus (FR) \sum_{p=1}^{m'} \int_{a_{p-1}(s)}^{a_p(s)} \int_{b_{p-1}(t)}^{b_p(t)}  k_p (s,t,\sigma,\tau) \odot \phi(z(\sigma,\tau)) d\tau d\sigma,
\end{equation}
and
\begin{equation}\label{2}
z(s,t) = g(s,t) \oplus (FR) \sum_{p=1}^{m'} \int_{a_{p-1}(s)}^{a_p(s)} \int_{b_{p-1}(t)}^{b_p(t)}  k_p (s,t,\sigma,\tau) \odot z(\sigma,\tau) d\tau d\sigma,
\end{equation}
where
$$
0= a_1 = a_0(s) < a_1(s) < ... < a_{m'-1}(s) < a_{m'}(s) = s,
$$
and
$$
0= a_1 = b_0(t) < b_1(t) < ... < b_{m'-1}(t) < b_{m'}(t) = t,
$$
the kernel $k_p: \Omega \times \Omega  \rightarrow R^+; p= 1, ..., m'$ are arbitrary positive functions, $s,t \in \Omega: [a_1,a_2] \times [a_1,a_2]$, $z,g: \Omega \rightarrow R_F$ are fuzzy functions, $\phi: R_F \rightarrow R_F$ is a continuous fuzzy function and $R_F$ shows the set of all fuzzy numbers.

Numerous researches have been done on fuzzy issues. In \cite{s2,s3}, Fariborzi Araghi and Noeiaghdam used numerical approaches to solve fuzzy integrals. By using the homotopy analysis transform technique and the Sinc-collocation approach, Noeiaghdam et al. have solved several types of fuzzy integral equations in \cite{s1,s4}. Noeiaghdam and Fariborzi Araghi have also explored the best outcomes and approach iterations for addressing fuzzy issues in \cite{s7,s8}. The fuzzy Taylor expansion was used by Allahviranloo et al. to solve fuzzy differential equations in \cite{s6} and Mikaeilvand et al. discussed the fuzzy system of equations in \cite{s5}.

The main goal of this work is to use the discrete collocation method to solve piecewise fuzzy Volterra integral equations which contain both linear and nonlinear problems. To discover the approximate solution to the issue, we use the two-dimensional fuzzy Lagrange interpolation polynomials and the two-dimensional fuzzy Guass-Legendre integration algorithm. By demonstrating the fundamental theorem and the corollary, we also talk about the solution's existence and uniqueness. Numerous instances in one and two-dimensional scenarios, both linear and nonlinear, are solved. We will demonstrate the accuracy and effectiveness of the strategy by plotting the graphs and tabulating the tables.

The following is how the paper has been set up: The two-dimensional fuzzy Gauss-Legendre quadrature formula and the two-dimensional fuzzy Lagrange interpolation are presented in Section 2. The appropriate discrete algebraic equations are found in Section 3 along with the discrete collocation method. The existence and uniqueness of the original equation's solution are established in Section 4. 

\section{Preliminaries}
The two-dimensional fuzzy Gauss-Legendre quadrature formula and the two-dimensional fuzzy Lagrange interpolation, which will be employed in this article, are introduced in this section.

\subsection{2-D fuzzy Lagrange interpolation polynomials}

Let $u(x, y) \in C_F (A =[a, b] \times [c, d]), C_F(A)$  represent the two-dimensional space of fuzzy continuous functions. To split $[a, b] \times [c, d]$ into $(N+1) \times (N+1)$ nodes, consider using $a =x_0< x_1< ... < x_N=b $ and $c=y_0<y_1<...<y_N=d$. Following that, $(x_i, y_j), 0 \leq i, j \leq N, N \in \mathbb{N}^*$ are the $(N+1)^2$ tensor product interpolation nodes on area $A =[a, b] \times [c, d]$. Assuming $I_Nu \in P^F_N$ is satisfied, and $P^F_N(A)$ is the space of all fuzzy algebraic polynomials of degree $N$ in each variable, $I_N$ will serve as the interpolation operator.

$$
(I_N u) (x_i,y_i) = u (x_i,y_j), ~~~ 0 \leq i,j \leq N
$$

The form of the two-dimensional fuzzy Lagrange interpolation polynomial is
\begin{equation}\label{3}
(I_N u) (x,y) = \sum_{i=0}^N  \sum_{j=0}^N  L_i (x)  L_j(y) \odot u(x_i,y_j)
\end{equation}
which $L_i(x)$ and $L_j(y)$ are the basis functions used for Lagrange interpolation, and they are defined as follows:

\begin{equation}\label{4}
L_i (x) = \prod_{k=0, k\neq i}^{N} \frac{(x-x_k)}{ (x_i - x_k)},~~~L_j (y) = \prod_{l=0, l\neq j}^{N} \frac{(y-y_l)}{ (y_j - y_l)},~~~ 0 \leq i,j \leq N. 
\end{equation}

In this study, we examine the possibility of employing the first class of Chebyshev nodes as interpolation nodes, namely
\begin{equation}\label{5}
x_i = \cos (\frac{2i+1}{2N+2} \pi),~~~ y_i = \cos (\frac{2j+1}{2M+2} \pi),~~~ 0 \leq i,j \leq N.
\end{equation}

\subsection{2-D Fuzzy Gauss-Legendre integration formula}
The formula for the fuzzy two-dimensional Gauss-Legendre quadrature is provided by

\begin{equation}\label{6}
(FR) \int_{-1}^{1} u(s,t) dt ds \approx \sum_{i=0}^{N}  \sum_{j=0}^{N} \omega_i \omega_j \odot u(s_i,t_j), 
\end{equation}
where the zeros of Legendre polynomials $L_{N+1}(x)$ are shown by $\{s_i\}^N_{i=0}$ and $\{t_j\}^N_{j=0}$  and 
\begin{equation}\label{7}
\omega_i = \frac{2}{(1-s_i^2) [L'_{N+1}(s_i)]^2},~~~\omega_j = \frac{2}{(1-t_j^2) [L'_{N+1}(t_j)]^2},~~~ 0 \leq i,j \leq N,
\end{equation}
are the corresponding weights.

\section{Main Method}

We use the following transformation
$$
s= \frac{1+x}{2}, t= \frac{1+y}{2}, x,y \in [-1,1],
$$
and
$$
\tau = \frac{1+\eta}{2}, \sigma = \frac{1+\xi}{2},  \xi \in [-1,x], \eta \in [-1,y]
$$
to rewrite Eqs. (\ref{1}) and (\ref{2}) as
\begin{equation}\label{19}
u(x,y) = f(x,y)  \oplus (FR) \sum_{p=1}^{m'} \int_{-1}^{a_p(x)} \int_{-1}^{b_p(y)}  \hat{k}_p (x,y,\xi,\eta) \odot H(u(\xi,\eta)) d\eta d\xi,
\end{equation}
and
\begin{equation}\label{20}
u(x,y) = f(x,y) \oplus (FR) \sum_{p=1}^{m'} \int_{-1}^{a_p(x)} \int_{-1}^{b_p(y)}  \hat{k}_p (x,y,\xi,\eta) \odot u(\xi,\eta) d\eta d\xi,
\end{equation}
where
$$
u(x,y) = z\left(\frac{1+x}{2}, \frac{1+y}{2}\right), ~~~ f(x,y) = g\left(\frac{1+x}{2}, \frac{1+y}{2}\right),
$$
$$
\hat{k}_p (x,y,\xi,\eta) = \left(\frac{1}{2}\right)^2  k_p \left(\frac{1+x}{2}, \frac{1+y}{2}, \frac{1+\xi}{2}, \frac{1+\eta}{2} \right),
$$
$$
H(u(\xi,\eta)) = \phi\left(z(\frac{1+\xi}{2}, \frac{1+\eta}{2})\right).
$$

For solving the nonlinear case (\ref{19}), if we substitute the interpolation points $(x_i,y_j)$ in Eq. (\ref{19}) we get
\begin{equation}\label{21}
u(x_i,y_j) = f(x_i,y_j)  \oplus (FR) \sum_{p=1}^{m'} \int_{-1}^{a_p(x_i)} \int_{-1}^{b_p(y_j)}  \hat{k}_p (x_i,y_j,\xi,\eta) \odot H(u(\xi,\eta)) d\eta d\xi,~~~0 \leq i,j \leq N.
\end{equation}

When we use the fuzzy Gauss-Legendre numerical integration, we should transform the integration domain $[-1, a_p(x_i)] \times [-1, b_p(y_j)]$ to $[-1,1]\times [-1,1]$ using the following transformation:
$$
\begin{array}{l}
  \xi = \xi(x_i, \theta) = \frac{1+a_p(x_i)}{2} \theta + \frac{a_p(x_i)-1}{2}  \\
  \\
  \eta = \eta(y_j, \theta) = \frac{1+b_p(y_j)}{2} \rho + \frac{b_p(y_j)-1}{2},~~~p=1,2,...,m',
\end{array}
$$
and Eq. (\ref{21}) can be written as
\begin{equation}\label{22}
u(x_i,y_j) = f(x_i,y_j)  \oplus (FR) \sum_{p=1}^{m'} \int_{-1}^{1} \int_{-1}^{1}  \tilde{k}_p (x_i,y_j,\xi(x_i,\theta),\eta(y_j,\rho)) \odot H(u(\xi(x_i,\theta),\eta(y_j,\rho))) d\rho d\theta,
\end{equation}
where
$$
\tilde{k}_p (x,y,\xi,\eta) = \frac{1+x}{2} . \frac{1+y}{2} \hat{k}_p (x,y,\xi,\eta).
$$

Applying the fuzzy Gauss-Legendre integration rule we get
\begin{equation}\label{23}
\tilde{u}(x_i,y_j) = f(x_i,y_j)  \oplus  \sum_{p=1}^{m'} \sum_{k=0}^{N} \sum_{l=0}^{N}  w_k w_l \odot  \tilde{k}_p (x_i,y_j,\xi(x_i,\theta_k),\eta(y_j,\rho_l)) \odot H(\tilde{u}(\xi(x_i,\theta_k),\eta(y_j,\rho_l))),
\end{equation}
where $\{ \theta_k \}_{k=0}^N $ and $\{ \rho_l \}_{l=0}^N$ are the sets of $N+1$ Gauss-Legendre quadrature points. Now we use the fuzzy Gauss-Legendre polynomial $U_n$ to interpolate the value of $\tilde{u}$
\begin{equation}\label{24}
U_N(x,y)  = \sum_{m=0}^{N} \sum_{n=0}^{N} L_m(x) L_n(y) \odot \tilde{u}(x_m,y_n).
\end{equation}

Combining Eqs. (\ref{23}) and (\ref{24}) we get
\begin{equation}\label{25}
\tilde{u}(x_i,y_j) = f(x_i,y_j)  \oplus  \sum_{p=1}^{m'} \sum_{k=0}^{N} \sum_{l=0}^{N}  w_k w_l \odot  \tilde{k}_p (x_i,y_j,\xi(x_i,\theta_k),\eta(y_j,\rho_l)) \odot H( \sum_{m=0}^{N} \sum_{n=0}^{N} L_m(\xi(x_i,\theta_k)) L_n(\eta(y_j,\rho_l)) \odot \tilde{u}(x_m,y_n) ),~~~ 0 \leq i,j  \leq N.
\end{equation}
If we solve the system of linear Eqs. (\ref{28}) then the values of  $\tilde{u}(x_i,y_j)$ can be found and substituting in $U_N (x,y) = \sum_{i=0}^{N} \sum_{j=0}^{N} L_i(x) L_j(y) \odot \tilde{u}(x_i,y_j)$ we will be able to find the approximate solution.

Repeating the process for the linear case (\ref{20}) and the interpolation points $(x_i,y_j)$ we can write
\begin{equation}\label{26}
u(x_i,y_j) = f(x_i,y_j) \oplus (FR) \sum_{p=1}^{m'} \int_{-1}^{a_p(x_i)} \int_{-1}^{b_p(y_j)}  \hat{k}_p (x_i,y_j,\xi,\eta) \odot u(\xi,\eta) d\eta d\xi,~~~0 \leq i,j \leq N.
\end{equation}
Applying the fuzzy Gauss-Legendre integration for the integral part of the above problem we get
\begin{equation}\label{26}
\tilde{u}(x_i,y_j) = f(x_i,y_j) \oplus  \sum_{p=1}^{m'} \sum_{k=0}^{N} \sum_{l=0}^{N}  w_k w_l \odot  \hat{k}_p (x_i,y_j,\xi(x_i,\theta_k),\eta(y_j,\rho_l)) \odot \tilde{u}(\xi(x_i,\theta_k),\eta(y_j,\rho_l)),
\end{equation}
and using the fuzzy Lagrange polynomials $U_n$ we have
\begin{equation}\label{28}
\tilde{u}(x_i,y_j) = f(x_i,y_j) \oplus  \sum_{m=0}^{N} \sum_{n=0}^{N} a_{m,n}  \odot  \tilde{u}(x_m,y_n),~~~0 \leq i,j \leq N.
\end{equation}
where
$$
a_{m,n} = \sum_{p=1}^{m'} \sum_{k=0}^{N} \sum_{l=0}^{N}  w_k w_l \odot  \hat{k}_p (x_i,y_j,\xi(x_i,\theta_k),\eta(y_j,\rho_l)) \odot L_m(\xi(x_i,\theta_k)) L_n(\eta(y_j,\rho_l)).
$$

Solving the system of linear Eqs. (\ref{28}) we can find the values of  $\tilde{u}(x_i,y_j)$ and the approximate solution can be obtained by $U_N (x,y) = \sum_{i=0}^{N} \sum_{j=0}^{N} L_i(x) L_j(y) \odot \tilde{u}(x_i,y_j)$.

\section{Existence and Uniqueness of the solution}

Assume that we have the following Lipschitz condition for the function $\phi(z(\sigma,\tau))$ as
\begin{equation}\label{8}
\exists L > 0; D (\phi(z(\sigma,\tau)), \phi(v(\sigma,\tau))) \leq L D (z(\sigma,\tau), v(\sigma,\tau)); \forall (\sigma,\tau) \in \Omega, z,v \in R_F,
\end{equation}
where $D (z(\sigma,\tau), v(\sigma,\tau))$ shows the Hausdorff distance between $z,v$.

\begin{lemma}\label{L1}
For nonnegative integrable functions  $k(s, t)$ and $y(s, t)$ the two-dimensional Gronwall inequality is defined as
\begin{equation}\label{9}
y(s,t) \leq C + \int_{0}^{s} \int_{0}^{t}  k (\sigma,\tau) y (\sigma,\tau)  d\tau d\sigma, ~~~\forall s,t \in \Omega
\end{equation}
and
\begin{equation}\label{10}
y(s,t) \leq C. exp ( \int_{0}^{s} \int_{0}^{t} k (\sigma,\tau) d\tau d\sigma ),
\end{equation}
where $C$ is a nonnegative constant.
\end{lemma}

Applying the condition (\ref{8}) and Lemma \ref{L1} we can prove the following theorem.

\begin{theorem}\label{T1}
Assume that  $k_p(s,t,\sigma,\tau), p=1,...,m'$ is a positive and continuous on $\Omega  \times \Omega$ and $\phi(z(\sigma,\tau)$ is continuous on $R_F$ which satisfies the condition (\ref{8}), then FVIE (\ref{1}) has a unique solution.
\end{theorem}

\textbf{Proof: } Assume that $z(s,t)$ and $v(s,t)$ are two different fuzzy solutions of (\ref{1}). Since $k_p(s,t,\sigma,\tau), p=1,...,m'$ sre continuous and positive functions on $\Omega  \times \Omega$ there exists $M_k^{(p)}>0$ such that $M_k^{(p)} = \max_{a_1 \leq s,t,\sigma,\tau \leq a_2} k_p(s,t,\sigma,\tau)$ and we get
\begin{equation}\label{11}
\begin{array}{l}
\displaystyle D(z(s,t), v(s,t)) \\
\\
\displaystyle = D \bigg(  g(s,t) \oplus (FR) \sum_{p=1}^{m'} \int_{a_{p-1}(s)}^{a_p(s)} \int_{b_{p-1}(t)}^{b_p(t)}  k_p (s,t,\sigma,\tau) \odot \phi(z(\sigma,\tau)) d\tau d\sigma, \\
\\
\displaystyle ~~~~~~~~~~g(s,t) \oplus (FR) \sum_{p=1}^{m'} \int_{a_{p-1}(s)}^{a_p(s)} \int_{b_{p-1}(t)}^{b_p(t)}  k_p (s,t,\sigma,\tau) \odot \phi(v(\sigma,\tau)) d\tau d\sigma   \bigg)\\
 \\
\displaystyle \leq \sum_{p=1}^{m'}   M_k^{(p)}  D \bigg(  (FR)   \int_{a_{p-1}(s)}^{a_p(s)} \int_{b_{p-1}(t)}^{b_p(t)} \phi(z(\sigma,\tau)) d\tau d\sigma, (FR)   \int_{a_{p-1}(s)}^{a_p(s)} \int_{b_{p-1}(t)}^{b_p(t)} \phi(v(\sigma,\tau)) d\tau d\sigma  \bigg)\\
 \\
\displaystyle \leq \sum_{p=1}^{m'}   M_k^{(p)} \int_{a_{p-1}(s)}^{a_p(s)} \int_{b_{p-1}(t)}^{b_p(t)} D (\phi(z(\sigma,\tau)), \phi(v(\sigma,\tau))) d\tau d\sigma\\
 \\
\displaystyle \leq \sum_{p=1}^{m'}   M_k^{(p)} L  \int_{a_{p-1}(s)}^{a_p(s)} \int_{b_{p-1}(t)}^{b_p(t)} D (z(\sigma,\tau), v(\sigma,\tau)) d\tau d\sigma, ~~~ \forall s,t \in \Omega
\end{array}
\end{equation}

We define the function $F: \Omega \rightarrow R^+ \cup \{ 0 \}$ and let $F(s,t) = D (z(s,t), v(s,t))$. From (\ref{11}) and Lemma \ref{L1} we get $F(s,t) = D (z(s,t), v(s,t))=0 $ and we have unique solution of problem (\ref{1}).

To show the existence of solution we use the following iterative process
\begin{equation}\label{13}
\begin{array}{l}
z_0(s,t) = g(s,t), \\
\\
z_n (s,t) = g(s,t) \oplus  (FR) \sum_{p=1}^{m'} \int_{a_{p-1}(s)}^{a_p(s)} \int_{b_{p-1}(t)}^{b_p(t)}  k_p (s,t,\sigma,\tau) \odot \phi(z_{n-1}(\sigma,\tau)) d\tau.
\end{array}
\end{equation}

Now we can prove the convergence of  the sequence $\{ z_n \}_{n=0}^{\infty}$.
\begin{equation}\label{14}
\begin{array}{l}
\displaystyle D(z_1(x,y), z_2 (x,y)) \\
\\
\displaystyle =D \bigg[ g(s,t) \oplus (FR) \sum_{p=1}^{m'} \int_{a_{p-1}(s)}^{a_p(s)} \int_{b_{p-1}(t)}^{b_p(t)}  k_p (s,t,\sigma,\tau) \odot \phi(z_{0}(\sigma,\tau)) d\tau d\sigma, \\
\\
\displaystyle g(s,t) \oplus (FR) \sum_{p=1}^{m'} \int_{a_{p-1}(s)}^{a_p(s)} \int_{b_{p-1}(t)}^{b_p(t)}  k_p (s,t,\sigma,\tau) \odot \phi(z_{1}(\sigma,\tau)) d\tau d\sigma   \bigg]\\
\\
\displaystyle \leq  \sum_{p=1}^{m'}  M_k^{(p)}  D \bigg[ (FR) \int_{a_{p-1}(s)}^{a_p(s)} \int_{b_{p-1}(t)}^{b_p(t)}  \phi(z_{0}(\sigma,\tau)) d\tau d\sigma,
\displaystyle (FR) \int_{a_{p-1}(s)}^{a_p(s)} \int_{b_{p-1}(t)}^{b_p(t)}  \phi(z_{1}(\sigma,\tau)) d\tau  d\sigma  \bigg]\\
\\
\displaystyle \leq  \sum_{p=1}^{m'}  M_k^{(p)}  \int_{a_{p-1}(s)}^{a_p(s)} \int_{b_{p-1}(t)}^{b_p(t)}  D\left( \phi(z_{0}(\sigma,\tau)),  \phi(z_{1}(\sigma,\tau))\right) d\tau d\sigma   \\
\\
\displaystyle \leq  \sum_{p=1}^{m'}  M_k^{(p)} L  \int_{a_{p-1}(s)}^{a_p(s)} \int_{b_{p-1}(t)}^{b_p(t)}  D\left( z_{0}(\sigma,\tau),  z_{1}(\sigma,\tau)\right) d\tau d\sigma  \\
\\
\displaystyle \leq  \sum_{p=1}^{m'}  M_k^{(p)} L s t D^*(z_0,z_1)
\end{array}
\end{equation}
where $D^*(z_0,z_1) = \sup_{a_1 \leq \sigma,\tau \leq a_2} D(z_{0}(\sigma,\tau),  z_{1}(\sigma,\tau))$.

\begin{equation}\label{15}
D(z_k(s,t), z_{k+1}(s,t)) \leq \frac{ (\sum_{p=1}^{m'}  M_k^{(p)} L s t)^k }{k!k!} D^*(z_0, z_1)
\end{equation}

thus for $n=k+1$ we get
\begin{equation}\label{16}
\begin{array}{l}
\displaystyle D(z_{k+1}(s,t), z_{k+2} (s,t)) \\
\\
\displaystyle =D \bigg[ g(s,t) \oplus (FR) \sum_{p=1}^{m'} \int_{a_{p-1}(s)}^{a_p(s)} \int_{b_{p-1}(t)}^{b_p(t)}  k_p (s,t,\sigma,\tau) \odot \phi(z_{k}(\sigma,\tau)) d\tau d\sigma, \\
\\
\displaystyle g(s,t) \oplus (FR) \sum_{p=1}^{m'} \int_{a_{p-1}(s)}^{a_p(s)} \int_{b_{p-1}(t)}^{b_p(t)}  k_p (s,t,\sigma,\tau) \odot \phi(z_{k+1}(\sigma,\tau)) d\tau d\sigma  \bigg]\\
\\
\displaystyle \leq  \sum_{p=1}^{m'}  M_k^{(p)}  D \bigg[ (FR) \int_{a_{p-1}(s)}^{a_p(s)} \int_{b_{p-1}(t)}^{b_p(t)}  \phi(z_{k}(\sigma,\tau)) d\tau d\sigma,
\displaystyle (FR) \int_{a_{p-1}(s)}^{a_p(s)} \int_{b_{p-1}(t)}^{b_p(t)}  \phi(z_{k+1}(\sigma,\tau)) d\tau d\sigma  \bigg]\\
\\
\displaystyle \leq  \sum_{p=1}^{m'}  M_k^{(p)}  \int_{a_{p-1}(s)}^{a_p(s)} \int_{b_{p-1}(t)}^{b_p(t)}  D\left( \phi(z_{k}(\sigma,\tau)),  \phi(z_{k+1}(\sigma,\tau))\right) d\tau d\sigma  \\
\\
\displaystyle \leq  \sum_{p=1}^{m'}  M_k^{(p)}   \int_{a_{p-1}(s)}^{a_p(s)} \int_{b_{p-1}(t)}^{b_p(t)}  D\left( z_{k}(\sigma,\tau),  z_{k+1}(\sigma,\tau)\right) d\tau  d\sigma \\
\\
\displaystyle \leq  \sum_{p=1}^{m'}  M_k^{(p)}   \int_{a_{p-1}(s)}^{a_p(s)} \int_{b_{p-1}(t)}^{b_p(t)}  \frac{ (\sum_{p=1}^{m'}  M_k^{(p)} L  \sigma \tau )^k }{k!k!} D^*(z_0, z_1)   d\tau  d\sigma \\
\\
\displaystyle \leq   \frac{ (\sum_{p=1}^{m'}  M_k^{(p)} L  s t  )^{k+1} }{(k+1)! (k+1)!} D^*(z_0, z_1).
\end{array}
\end{equation}
Thus the inequality (\ref{15}) holds for all $n \in  N$. For any positive integer $q$ we get
\begin{equation}\label{17}
\begin{array}{l}
\displaystyle D(z_n(s,t), z_{n+q}(s,t) ) \leq  D(z_n(s,t), z_{n+1}(s,t) ) + ... +  D(z_{n+q-1}(s,t), z_{n+q}(s,t) )\\
 \\
\displaystyle \leq  \left(\frac{ (\sum_{p=1}^{m'}  M_k^{(p)} L  s t  )^{n} }{(n)! (n)!} + ... + \frac{ (\sum_{p=1}^{m'}  M_k^{(p)} L  s t  )^{n+q-1} }{(n+q-1)! (n+q-1)!} \right) D^*(z_0, z_1),~~~ \forall s,t \in \Omega.
\end{array}
\end{equation}

From (\ref{17}) we have:
\begin{equation}\label{18}
  \lim_{n \rightarrow 0} D(z_n(s,t), z_{n+q}) = 0.
\end{equation}
That is the fuzzy sequence $\{ z_n(s,t)  \}_{n=0}^{\infty}$ is convergent and its limit is the only solution of (\ref{1}).

\begin{corollary}
  Let  $k_p (s,t,\sigma,\tau), p= 1,...,m'$ be continuous and positive on $\Omega  \times \Omega$ and $z(s,t)$ be continuous on $\Omega$ the linear two dimensional VIE (\ref{2}) has a unique solution.
\end{corollary}

\textbf{Proof:}  Same as the prove of Theorem \ref{T1}.

\section{Conclusion}
This study has investigated the existence and uniqueness of solutions to Fuzzy Volterra integral equations (FVIEs) with discontinuous kernels, which have significant applications in load leveling problems in energy storage. Both linear and nonlinear cases were examined, and the collocation method was employed to analyze the equations. The research aimed to provide a deeper understanding of the solution properties of FVIEs with discontinuous kernels, which is essential for developing effective numerical methods. The findings of this study contribute to the advancement of mathematical modeling and simulation in energy storage systems, enabling more accurate predictions and optimized performance. The results have important implications for the design and operation of energy storage systems, and can inform decision-making in this field.




\begin{thebibliography}{99}


\bibitem{21} D. Juraev, S. Noeiaghdam, Modern Problems of Mathematical Physics and Their Applications, MDPI, March 2022, ISBN 978-3-0365-3496-1 (Hbk), ISBN 978-3-0365-3495-4 (PDF).

\bibitem{48} D. Sidorov, I. Muftahov, D. Karamov, N. Tomin, D. Panasetsky, A. Dreglea, F. Liu, A. Foley, A dynamic analysis of energy storage with renewable and diesel generation using Volterra equations, IEEE Trans. Ind. Inform. 16(5) (2020) 3451–3459.

\bibitem{63} Q. Tao, F. Liu, D. Sidorov, Recurrent neural networks application to forecasting with two cases: load and pollution, in: G.W. Weber (Ed.), Intelligent Systems and Computing, vol.1072, Springer, 2020, pp.369–378.

\bibitem{49} D. Sidorov, A. Tynda, I. Muftahov, A. Dreglea, F. Liu, Nonlinear systems of Volterra equations with piecewise smooth kernels: numerical solution and application for power systems operation, Mathematics 8(8) (2020) 1257.

\bibitem{12} A. Domyshev, D. Sidorov, D. Panasetsky, An improved two-stage optimization procedure for optimal power flow, Energy Syst. Res. 3(1) (2020) 9, 9pp.

\bibitem{57} D.N. Sidorov, A.V. Zhukov, I.R. Muftahov, Volterra equation based models for energy storage usage based on load forecast in EPS with renewable generation, Bull. Irkutsk State Univ., Ser. Math. 26 (2018) 76–90.

\bibitem{37} S. Noeiaghdam, D. Sidorov, I. Muftahov, A.V. Zhukov, Control of accuracy on Taylor-collocation method for load leveling problem, Bull. Irkutsk State Univ., Ser. Math. 30 (2019) 59–72, https://doi .org /10 .26516 /1997 -7670 .2019 .30 .59.



\bibitem{s1} M. A. Fariborzi Araghi, S. Noeiaghdam, Homotopy analysis transform method for solving generalized Abel's fuzzy integral equations of the first kind, 4th Iranian Joint Congress on Fuzzy and Intelligent Systems, CFIS 2015, 2016, 7391645. https://doi.org/10.1109/CFIS.2015.7391645

\bibitem{s2}M. A. Fariborzi Araghi, S. Noeiaghdam, Finding the optimal step of fuzzy Newton- Cotes integration rules by using CESTAC method, Journal of Fuzzy Set Valued Analysis, 2 (2017) 62-85.

\bibitem{s3} M. A. Fariborzi Araghi, S. Noeiaghdam, A valid scheme to evaluate fuzzy definite integrals by applying the CADNA library, International Journal of Fuzzy System Applications, 6 (4) (2017) 1-20.

\bibitem{s4}S. Noeiaghdam, M.A. Fariborzi Araghi, S. Abbasbandy, Valid implementation of Sinc-collocation method to solve the fuzzy Fredholm integral equation, Journal of Computational and Applied Mathematics, 370 (2020) 112632. https://doi.org/10.1016/j.cam.2019.112632

\bibitem{s5}N. Mikaeilvand, Z. Noeiaghdam, S. Noeiaghdam, Juan J. Nieto, A novel technique to solve the fuzzy system of equations, Mathematics, 2020, 8(5), 850.  https://doi.org/10.3390/math8050850

\bibitem{s6}T. Allahviranloo, Z. Noeiaghdam, S. Noeiaghdam, Juan J. Nieto , A Fuzzy Method for Solving Fuzzy Fractional Differential Equations Based on the Generalized Fuzzy Taylor Expansion, Mathematics, 2020, 8(12), 1–24, 2166.   https://doi.org/10.3390/math8122166

\bibitem{s7}M.A. Fariborzi Araghi, S. Noeiaghdam (2022) Finding Optimal Results in the Homotopy Analysis Method to Solve Fuzzy Integral Equations. In: Allahviranloo T., Salahshour S. (eds) Advances in Fuzzy Integral and Differential Equations. Studies in Fuzziness and Soft Computing, vol 412. Springer, Cham. (Scopus) https://doi.org/10.1007/978-3-030-73711-5$_{-}$7

\bibitem{s8} S. Noeiaghdam, M. A. Fariborzi Araghi, (2021) Application of the CESTAC Method to Find the Optimal Iteration of the Homotopy Analysis Method for Solving Fuzzy Integral Equations. In: Allahviranloo T., Salahshour S., Arica N. (eds) Progress in Intelligent Decision Science. IDS 2020. Advances in Intelligent Systems and Computing, vol 1301. Springer, Cham.  https://doi.org/10.1007/978-3-030-66501-2$_{-}$49

\bibitem{s100}S. Noeiaghdam, D. Sidorov, V. Sizikov, N. Sidorov, Control of accuracy on Taylor-collocation method to solve the weakly regular Volterra integral equations of the first kind by using the CESTAC method, Applied and Computational Mathematics an International Journal, 19 (1) (2020) 81-105. 

\bibitem{s101}S. Noeiaghdam, A. Dreglea, J. H. He, Z. Avazzadeh, M. Suleman, M. A. Fariborzi Araghi, D. Sidorov, N. Sidorov, Error estimation of the homotopy perturbation method to solve second kind Volterra integral equations with piecewise smooth kernels: Application of the CADNA library, Symmetry, 2020, 12(10), 1–16, 1730.   https://doi.org/10.3390/sym12101730 

\bibitem{s102}E. Hashemizadeh, M. A. Ebadi, S. Noeiaghdam, Matrix method by Genocchi polynomials for solving nonlinear Volterra integral equations with weakly singular kernel, Symmetry, 2020, 12(12), 1–19, 2105. https://dx.doi.org/10.3390/sym12122105

\bibitem{s103}S. Noeiaghdam, D. Sidorov, A. M. Wazwaz, N. Sidorov, V. Sizikov, The numerical validation of the Adomian decomposition method for solving Volterra integral equation with discontinuous kernel using the CESTAC method, Mathematics, 2021, 9(3), 1–15, 260. https://doi.org/10.3390/math9030260

\bibitem{s104}E. Zarei, S. Noeiaghdam, Advantages of the Discrete Stochastic Arithmetic to Validate the Results of the Taylor Expansion Method to Solve the Generalized Abel’s Integral Equation. Symmetry 2021, 13, 1370. https://doi.org/10.3390/sym13081370

\bibitem{s105}S. Noeiaghdam, S. Micula, A Novel Method for Solving Second Kind Volterra Integral Equations with Discontinuous Kernel. Mathematics 2021, 9, 2172. https://doi.org/10.3390/math9172172

\bibitem{s106}S. Noeiaghdam, M. A. Fariborzi Araghi, D. Sidorov, Dynamical strategy on homotopy perturbation method for solving second kind integral equations using the CESTAC method, Journal of Computational and Applied Mathematics, 2022, 114226. https://doi.org/10.1016/j.cam.2022.114226

\bibitem{s107}S. Noeiaghdam, D. Sidorov, A. Dreglea, A novel numerical optimality technique to find the optimal results of Volterra integral equation of the second kind with discontinuous kernel, Applied Numerical Mathematics 186 (2023) 202–212. https://doi.org/10.1016/j.apnum.2023.01.011

\end{thebibliography}
\end{document}